\documentclass[12pt]{amsart}
\usepackage{amsmath}
\usepackage{amssymb}
\usepackage{latexsym}
\usepackage{amscd}
\usepackage{citesort}

\newdimen\AAdi%
\newbox\AAbo%
\def\AAk#1#2{\s_etbox\AAbo=\hbox{#2}\AAdi=\wd\AAbo\kern#1\AAdi{}}%
\def\AAr#1#2#3{\s_etbox\AAbo=\hbox{#2}\AAdi=\ht\AAbo\raise#1\AAdi\hbox{#3}}%
\font\tenmsb=msbm10 at 12pt
\font\sevenmsb=msbm7 at 8pt
\font\fivemsb=msbm5 at 6pt
\newfam\msbfam
\textfont\msbfam=\tenmsb
\scriptfont\msbfam=\sevenmsb
\scriptscriptfont\msbfam=\fivemsb

\newcommand{\ba}{\begin{array}}
\newcommand{\ea}{\end{array}}

\newcommand{\Section}[2]{\setcounter{equation}{0}
\allowdisplaybreaks
\section[#1]{#2}}

\begin{document}
\title
[COMPARISON THEOREMS IN FINSLER GEOMETRY AND THEIR APPLICATIONS]
{COMPARISON THEOREMS IN FINSLER GEOMETRY AND THEIR APPLICATIONS}
\author
[B.Y. Wu]{B.Y. Wu}
\address[B.Y. Wu]
{Department of Mathematics\\ Zhejiang Normal University\\
Jinhua\\Zhejiang \\
and Institute of Mathematics\\  Fudan University\\  Shanghai }
\email{wubingye@zjnu.cn}

\author
[Y.L. Xin]{Y. L. Xin}
\address[Y.L. Xin]
{Institute of Mathematics\\  Fudan University\\  Shanghai}
\email{ylxin@fudan.edu.cn}
\thanks{The research of the second author was partially supported by
NSFC of and SFECC}

\renewcommand{\subjclassname}{%
  \textup{2000} Mathematics Subject Classification}
\subjclass{Primary 53C60; Secondary 53B40 }\date{}
\maketitle

\begin{abstract}
We prove Hessian comparison
theorems, Laplacian comparison theorems and volume comparison
theorems of Finsler manifolds under various curvature conditions.
As applications, we derive Mckean type theorems for the first eigenvalue of Finsler
manifolds, as well as generalize a result on fundamental group due
to  Milnor to Finsler manifolds.
\end{abstract}

\Section{Introduction }{Introduction }

Finsler geometry, a natural generalization of
Riemannian geometry, was initiated by Finsler P. [Fin] in 1918,
from considerations of regular problems in the calculus of
variations. It developed steadily, with much investigation from
the geometric point of view. Chern [Ch1] and many others defined
various connections in Finsler manifolds, along the lines of the
Levi-Civita connection in Riemannian manifolds; for a comprehensive account,
see [BCS].

Recently, there has been a surge of interest in
Finsler geometry, especially in its global and analytic aspects. A
natural question, that has lately attracted some attention (see
e.g., [AL,Ce,Sh2]), is how to generalize the Laplacian from
Riemannian manifolds to Finsler manifolds. In Riemannian case, the
Laplacian of a function equals to the divergence of gradient of
the function, and the spectrum of Laplacian on Riemannian
manifolds has been extensively studied. we shall adopt the notion of
Laplacian for Finsler manifolds used in [Sh2].

The comparison technique
is is widely used in Riemannian geometry. To pursue the global Finsler
geometry we would generalize comparison theorems to the Finsler
setting. It has been started in [Sh3]. The present paper would continue
the investigation on this direction. We derive
some Hessian comparison theorems, Laplacian comparison theorems
and volume comparison theorems of Finsler manifolds under the various
curvature assumptions.  Then, we
give some applications. We obtain some Mckean type theorems for
the first eigenvalue of Finsler manifolds, as well as generalize a
result on fundamental group due to Milnor to Finsler manifolds.

\Section{Finsler Geometry}{Finsler Geometry}
Let $(M,F)$ be a Finsler $n$-manifold with Finsler
metric $F:TM\rightarrow [0,\infty)$. Let $(x,y)=(x^i,y^i)$ be the
local coordinates on $TM$, and $\pi:TM\backslash 0\rightarrow M$
the natural projection. Unlike in the Riemannian case, most
Finsler quantities are functions of $TM$ rather than $M$. Some
frequently used quantities and relations:
$$g_{ij}(x,y):=\frac{1}{2}\frac{\partial^2F^2(x,y)}{\partial y^i\partial
y^j},\quad {\rm (fundamental\quad tensor)}$$
$$C_{ijk}(x,y):=\frac{1}{4}\frac{\partial^3F^2(x,y)}{\partial y^i\partial
y^j\partial y^k}, \quad {\rm (Cartan\quad tensor)}$$
$$(g^{ij}):=(g_{ij})^{-1},$$
$$\gamma^k_{ij}:=\frac{1}{2}g^{km}\left( \frac{\partial g_{mj}}{\partial
x^i} +\frac{\partial g_{im}}{\partial x^j}-\frac{\partial
g_{ij}}{\partial x^m}\right),$$
$$N^i_j=\gamma^i_{jk}y^k-C^i_{jk}\gamma^k_{rs}y^ry^s.$$
According to [Ch1], the pulled-back bundle $\pi^*TM$ admits a
unique linear connection, called the {\sl Chern connection}. Its
connection forms are characterized by the structure equation:\\
$\bullet$  Torsion freeness:
$$dx^j\wedge\omega^i_j=0;$$
$\bullet$ Almost $g$-compatibility:
$$dg_{ij}-g_{kj}\omega^k_i-g_{ik}\omega^k_j=2C_{ijk}(dy^k+N^k_ldx^l).$$
It is easy to know that torsion freeness is equivalent to the
absence of $dy^k$ terms in $\omega^i_j$; namely,
$$\omega^i_j=\Gamma^i_{jk}dx^k,$$
together with the symmetry
$$\Gamma^i_{jk}=\Gamma^i_{kj}.$$
Let $V=v^i\partial/\partial x^i$ be a non-vanishing vector field
on an open subset $\mathcal{U}\subset$ $M$. One can introduce a
Riemannian metric $g_V$ and a linear connection $\nabla^V$ on the
tangent bundle over $\mathcal{U}$ as following:
$$g_V(X,Y):=X^iY^jg_{ij}(x,v),\quad \forall X=X^i\frac{\partial}{\partial
x^i},Y=Y^i\frac{\partial}{\partial x^i};$$
$$\nabla^V_{\frac{\partial}{\partial x^i}}\frac{\partial}{\partial
x^j}:=\Gamma^k_{ij}(x,v)\frac{\partial}{\partial x^k}.$$ From the
torsion freeness and $g$-compatibility of Chern connection we have
$$\nabla^V_XY-\nabla^V_YX=[X,Y],\eqno{(2.1)}$$
$$Xg_V(Y,Z)=g_V(\nabla^V_XY,Z)+g_V(Y,\nabla^V_XZ)+2C_V(\nabla^V_XV,Y,Z), \eqno{(2.2)}$$
here $C_V$ is defined by
$$C_V(X,Y,Z)=X^iY^jZ^kC_{ijk}(x,v),$$ and it satisfies
$$C_V(V,X,Y)=0.\eqno{(2.3)}$$ The  {\sl Chern curvature}
$R^V(X,Y)Z$ for vector fields $X,Y,Z$ on $\mathcal{U}$ is defined
by
$$R^V(X,Y)Z:=\nabla^V_X\nabla^V_YZ-\nabla^V_Y\nabla^V_XZ-\nabla^V_{[X,Y]}Z.$$
In the Riemannian case this curvature does not depend on $V$ and
coincides with the Riemannian curvature tensor. For a {\sl flag}
$(V;\sigma)$ (or $(V;W)$) consisting of a non-zero tangent vector
$V\in T_xM$ and a 2-plane $\sigma\subset T_xM$ with $V\in \sigma$
the {\sl flag curvature} $K(V;\sigma)$ is defined as following:
$$K(V;\sigma)=K(V;W):=\frac{g_V(R^V(V,W)W,V)}{g_V(V,V)g_V(W,W)-g_V(V,W)^2}.$$
Here $W$ is a tangent vector, such that $V,W$ span the 2-plane
$\sigma$ and $V\in T_xM$ is extended to a geodesic field, i.e.,
$\nabla^V_VV=0$ near $x$. In the Riemannian case the flag
curvature is the sectional curvature of the 2-plane $\sigma$ and
does dot depend on $V$. In the literature there are several
connections used in Finsler geometry, but for the definition of
the flag curvature it does not make a difference whether one uses
the Chern, the Cartan or the Berwald connection. The {\sl Ricci
curvature} of $V$ is defined by
$$ Ric(V)=\sum_iK(V;E_i),$$ where $E_1,\cdots,E_n$ is the local
$g_V$-orthonormal frame over $\mathcal{U}$.\\
\hspace*{0.7cm} A Finsler metric $F$ on $M$ is called {\sl
reversible} if $F(-X)=F(X)$ for all $X\in TM$. In order to
consider the {\sl non-reversible} Finsler metric Rademacher [Ra]
introduced the {\sl reversibility} $\lambda=\lambda(M,F)$ as
following
$$\lambda:=\sup_{X\in TM\backslash 0}\frac{F(-X)}{F(X)}.\eqno{(2.4)}$$
Clearly $\lambda\in [1,\infty]$ and $\lambda=1$ if and only if $F$
is reversible.\\
\hspace*{0.7cm} Let $\gamma(t),0\le t\le l$ be a geodesic with
unit speed velocity field $T$. A vector field $J$ along $\gamma$
is called to be a {\sl Jacobi field} if it satisfies the following
equation
$$\nabla^T_T\nabla^T_TJ+R^T(J,T)T=0.$$
For vector fields $X$ and $Y$ along $\gamma$, the {\sl index form}
$I_{\gamma}(X,Y)$ is defined by
$$I_{\gamma}=\int_{0}^{l}\left(g_T(\nabla^T_TX,\nabla^T_TY)-g_T(R^T(X,T)T,Y)\right)dt.$$
 Let
$${\rm ct}_c(t)=\left\{ \begin{array}{ll}
\sqrt{c}\cdot {\rm cotan}(\sqrt{c}t), & c>0\\ \frac{1}{t}, & c=0\\
\sqrt{-c}\cdot {\rm cotanh}(\sqrt{-c}t), & c<0
\end{array}\right. \eqno{(2.5)}$$ The following result is
fundamental.\\
{\bf Lemma 2.1}([BCS],page 254) {\sl Let $(M,F)$ be an Finsler
manifold and $\gamma(t),0\le t\le l$ be a geodesic with unit speed
velocity
field $T(t)$. Suppose: }\\
$\bullet$ {\sl The flag curvature $K(T;W)\le c$ for any $W\in
T_{\gamma(t)}M$.}\\
$\bullet$ {\sl $J$ is a Jacobi field along $\gamma$ that is
$g_T$-orthogonal to $\gamma$. }\\
$\bullet$ $J(0)=0.$\\
{\sl Then for $0<t\le l$ when $c\le 0$ or
$0<t<\frac{\pi}{\sqrt{c}}$ when $c>0$, }
$$\left.\frac{g_T(\nabla^T_TJ,J)}{g_T(J,J)}\right|_t\ge
{\rm ct}_c(t).$$

\Section{Laplacian for Finsler Manifolds}{Laplacian for Finsler Manifolds}
In this section we shall introduce the Laplacian
for Finsler manifolds adopted in [Sh2]. For this purpose, let us
first recall the notion of Legendre transformation.\\
\hspace*{0.7cm} Given a Finsler manifold $(M,F)$, the {\sl dual
Finsler metric} $F^*$ on $M$ is defined by
$$F^*(\xi_x):=\sup_{Y\in T_xM\backslash
0}\frac{\xi(Y)}{F(Y)},\quad\forall \xi\in T^*M,$$ and the
corresponding fundamental tensor is defined by
$$g^{*kl}(\xi)=\frac{1}{2}\frac{\partial^2F^{*2}(\xi)}{\partial\xi_k\partial\xi_l}.$$
The {\sl Legendre transformation} $l:TM\rightarrow T^*M$ is
defined by
$$l(Y)=\left\{ \begin{array}{ll}
g_Y(Y,\cdot), & Y\ne 0\\ 0, & Y=0.
\end{array}\right.$$
The following result is well-known (see [BCS, Sh2]).\\
{\bf Lemma 3.1} {\sl For any $x\in M$, the Legendre transformation
is a smooth diffeomorphism from $T_xM\backslash 0$ onto
$T^*_xM\backslash 0$, and it is norm-preserving, namely,
$F(Y)=F^*(l(Y)),\forall Y\in TM$. Consequently,} $g^{ij}(Y)=g^{*ij}(l(Y)).$\\
\hspace*{0.7cm} Now let $f:M\rightarrow \mathbb{R}$ be a smooth
function on $M$. The {\sl gradient} of $f$ is defined by $\nabla
f=l^{-1}(df)$. Thus we have
$$ df(X)=g_{\nabla f}(\nabla f,X),\quad X\in TM.$$
Let $\mathcal{U}=$ $\{ x\in M:\nabla f\mid_x\ne 0\}$. We define
the {\sl Hessian} $H(f)$ of $f$ on $\mathcal{U}$ as following:
$$H(f)(X,Y):=XY(f)-\nabla^{\nabla f}_XY(f),\quad \forall X,Y\in
TM\mid_{\mathcal{U}}.\eqno{(3.1)}$$ From (2.1)-(2.3) we see that
$H(f)$ is symmetric, and it can be rewritten as
$$ H(f)(X,Y)=g_{\nabla f}(\nabla^{\nabla f}_X\nabla
f,Y).\eqno{(3.2)}$$ It should be noted here that the notion of
Hessian here is different from that in [Sh1-2]. In that case
$H(f)$ is in fact defined by
$$H(f)(X,X)=XX(f)-\nabla^X_XX(f),$$ and there is no definition for
$H(f)(X,Y)$ if $X\ne Y$.\\
\hspace*{0.7cm} In order to define the divergence for vector
field, we need the volume form on $M$. A {\sl volume form} $d\mu$
on $M$ is nothing but a global non-degenerate $n$-form on $M$. A
frequently used volume form for $(M,F)$ is the so-called {\sl
Busemann-Hausdorff volume form} $dV_F$ which is locally expressed
by $dV_F=\sigma_F(x)dx^1\wedge\cdots\wedge dx^n$, where
$$\sigma_F(x):=\frac{{\rm vol}(B^n(1))}{{\rm vol}\left((y^i)\in R^n: F(x,y^i\frac{\partial}{\partial x^i})
<1\right)}.$$
 In the
following we consider the Finsler manifold $(M,F,d\mu)$ equipped
with a volume form $d\mu$. Let $X\in TM$. The {\sl divergence}
div$(X)$ of $X$ is defined by
$$d(X\rfloor d\mu)={\rm div}(X)d\mu.\eqno{(3.3)}$$ In local
coordinate system $(x^i)$, express
$d\mu=\sigma(x)dx^1\wedge\cdots\wedge dx^n$. Then for vector field
$X=X^i\partial/\partial x^i$ on $M$,
$${\rm div}(X)=\frac{1}{\sigma}\frac{\partial}{\partial x^i}\left(
\sigma X^i\right).\eqno{(3.4)}$$ Applying the Stokes theorem to
$\eta=X\rfloor d\mu$ we have\\
{\bf Lemma 3.2} ([Sh1-2]) {\sl Let $(M,F,d\mu)$ be a Finsler
$n$-manifold. Let $\Omega$ be a compact domain with smooth
boundary $\partial\Omega$ and $\nu$ denote the outward pointing
normal vector. Then for any smooth vector field $X$ on $M$, }
$$\int_{\Omega}{\rm
div}(X)d\mu=\int_{\partial\Omega} g_{\nu}(\nu,X)dA_{\mu},$$ {\sl
where $dA_{\mu}$ is the volume form on $\partial\Omega$ induced
from $d\mu$.}\\
For $y\in T_xM\backslash 0$, define
$$\tau(y):=\log\frac{\sqrt{\det\left(g_{ij}(x,y)\right)}}{\sigma}.\eqno{(3.5)}$$
$\tau$ is called the {\sl distorsion} of $(M,F,d\mu)$. To measure
the rate of the distorsion along geodesics, we define
$${\bf
S}(y):=\frac{d}{dt}\left[\tau(\dot{\gamma}(t))\right]_{t=0},\eqno{(3.6)}$$
where $\gamma(t)$ is the geodesic with $\dot{\gamma}(0)=y$. {\bf
S} is called the {\sl S-curvature}[Sh2], and it is an important
non-Riemannian curvature for Finsler manifold. In local
coordinates it can be expressed by[Sh2]
$${\bf S}(y)=N^i_i(x,y)-\frac{y^i}{\sigma(x)}\frac{\partial\sigma}{\partial
x^i}(x).\eqno{(3.7)}$$ Now we are ready to introduce the Laplacian
$\triangle f$ of $f$ as $\triangle f=$div$(\nabla
f)=$div$\left(l^{-1}(df)\right)$. By Lemma 3.1 and (3.4) we have
the following local express for $\triangle f$.
$$\triangle f=\frac{1}{\sigma(x)}\frac{\partial}{\partial x^i}\left(
\sigma(x)g^{*ij}(df)\frac{\partial f}{\partial x^j}\right)$$
$$=\frac{1}{\sigma(x)}\frac{\partial}{\partial x^i}\left(
\sigma(x)g^{ij}(\nabla f)\frac{\partial f}{\partial
x^j}\right).\eqno{(3.8)}$$ For later use we need the following
invariant express for $\triangle f$.\\
{\bf Lemma 3.3} {\sl Let $(M,F,d\mu)$ be a Finsler $n$-manifold,
and $f:M\rightarrow\mathbb{R}$ the smooth function on $M$. Then on
$\mathcal{U}=$ $\{ x\in M:\nabla f\mid_x\ne 0\}$ we have}
$$\triangle f=\sum_a H(f)(e_a,e_a)-{\bf S}(\nabla f):={\rm tr}_{\nabla
f}H(f)-{\bf S}(\nabla f),$$ {\sl where $e_1,\cdots,e_n$ is the
local $g_{\nabla f}$-orthonormal frame on} $\mathcal{U}$.\\
{\bf Proof.} Write
$$e_a=u_a^i\frac{\partial}{\partial x^i},\quad \frac{\partial}{\partial
x^i}=v^a_ie_a,\eqno{(3.9)}$$ then
$$v^a_iu^i_b=\delta^a_b,\quad u_a^iv_j^a=\delta^i_j,\quad
g^{ij}=\sum_au_a^iu_a^j.\eqno{(3.10)}$$ substituting (3.9) and
(3.10) into (3.8) we get
$$\Delta f=\sum_ae_ae_a(f)+(\nabla f)(\log\sigma)+\sum_be_a(u_b^i)
v_i^ae_b(f).\eqno{(3.11)}$$ From (3.7) one has
$${\bf S}(\nabla f)=N^i_i(\nabla f)-(\nabla
f)(\log\sigma).\eqno{(3.12)}$$ Let $\{\omega^b_a\}$ be the Chern
connection form with respect to $\{e_a\}$, then it is easy to
deduce that (see [BCS], page 42)
$$\omega^a_b=(du_b^i)v^a_i+u_b^j\omega^i_jv_i^a,\eqno{(3.13)}$$
$$\omega^a_b+\omega^b_a=-2C_{abc}v^c_i(dy^i+N^i_jdx^j).\eqno{(3.14)}$$
Thus from (3.13) and (3.14) we have
$$\sum_a\nabla^{\nabla f}_{e_a}e_a=\sum_a\omega_a^b(e_a)e_b$$
$$=\sum_b\left(-e_a(u_b^i)v_i^a-u_b^j\omega_j^i\left(\frac{\partial}{\partial
x^i}\right)
-2\sum_aC_{abc}v^c_iN^i_ju^j_a\right)e_b.\eqno{(3.15)}$$ Noting
that $\Gamma^i_{jk}(\nabla f)^k=N^i_j(\nabla f)$ ([BCS], page43),
we deduce from (2.3),(3.9),(3.10) and (3.15) that
$$\sum_a\nabla^{\nabla
f}_{e_a}e_a(f)=-\sum_b\left(e_a(u_b^i)v_i^ae_b(f)+\Gamma^i_{ji}u_b^ju_b^k\frac{\partial
f}{\partial x^k}\right)$$
$$=-\sum_be_a(u_b^i)v_i^ae_b(f)-N^i_i(\nabla f).\eqno{(3.16)}$$
Combining (3.1),(3.11),(3.12) and (3.16) we obtain the desired
result.\\

\Section{The Hessian Comparison Theorem}{The Hessian Comparison Theorem}
In this section let us study the Hessian
comparison theorem for distance function. For this purpose, let
us first compute the Hessian of distance function.\\
\hspace*{0.7cm} Let $(M,F,d\mu)$ be a Finsler $n$-manifold, and
$r=d_F(p,\cdot)$ is the distance function on $M$ from a fixed
point $p\in M$. It is well-known that $r$ is smooth on
$M\backslash \{p\}$ away from the cut points of $p$. Now we assume
$\gamma$ is a unit-speed geodesic without a conjugate point up to
distance $r$ from $p$. It is known that $F(\nabla r)=1$ (see
[Sh2],page 38), which together with the first variation of arc
length implies that  $\nabla r=T:=\dot{\gamma}$. For any vector
$X\in T_{\gamma(r)}M$, there exists a unique Jacobi field $J$ such
that $J(0)=0,J(r)=X$. We have, by (2.1)-(2.3) and (3.2),
$$H(r)(X,X)=g_T(\nabla^T_XT,X)=g_T(\nabla^T_JT,J)\mid^{\gamma(r)}_{\gamma(0)}=\int_0^r\frac{d}{dt}
g_T(\nabla^T_JT,J)dt$$
$$=\int_0^r\left(g_T(\nabla^T_T\nabla^T_JT,J)+g_T(\nabla^T_JT,\nabla^T_TJ)\right)dt$$
$$=\int_0^r\left(g_T(R^T(T,J)T,J)+g_T(\nabla^T_TJ,\nabla^T_TJ)\right)dt\eqno{(4.1)}$$
$$=I_{\gamma}(J,J)=I_{\gamma}(J^{\bot},J^{\bot})=g_T(\nabla^T_TJ^{\bot},J^{\bot}),$$
where $J^{\bot}=J-g_T(T,J)T$. Now we can prove the following
Hessian comparison theorem.\\
{\bf Theorem 4.1} {\sl Let $(M,F,d\mu)$ be a Finsler $n$-manifold,
$r=d_F(p,\cdot)$, the distance function from a fixed point $p$.
Suppose that the flag curvature of $M$ satisfies $K(V;W)\le c$
(resp. $K(V;W)\ge c$) for any $V,W\in TM$. Then for any vector $X$
on $M$ the following inequality holds whenever $r$ is smooth:}
$$H(r)(X,X)\ge({\sl resp.}\le){\rm ct}_c(r)\left(g_{\nabla r}(X,X)-g_{\nabla
r}(\nabla r,X)^2\right).$$ {\bf Proof.} First we note that
$$g_T(J^{\bot},J^{\bot})\mid_{\gamma(r)}=g_T(X,X)-g_T(T,X)^2,$$
by (4.1) and Lemma 2.1 we conclude that in the case $K(V;W)\le c$
one has
$$H(r)(X,X)\ge{\rm ct}_c(r)\left(g_{\nabla r}(X,X)-g_{\nabla
r}(\nabla r,X)^2\right).$$ Now we consider the case $K(V;W)\ge c$.
For given $X$, by parallel transformation along $\gamma$ we obtain
a vector field $X(t)$ along $\gamma$. We define a vector field
$W(t)$ along $\gamma$ by  $W(t)=\frac{{\rm s}_c(t)}{{\rm
s}_c(r))}X(t)$, where
$${\rm s}_c(t)=\left\{ \begin{array}{ll}
\sin(\sqrt{c}t), & c>0\\ t, & c=0\\
\sinh(\sqrt{-c}t), & c<0
\end{array}\right.  .\eqno{(4.2)}$$
It is clear that $W(0)=J(0)=0,W(r)=J(r)$, and consequently,
$W^{\bot}(0)=J^{\bot}(0)=0,W^{\bot}(r)=J^{\bot}(r)$ . Thus from
(4.1) and the basic index lemma (see [BCS],page 182) we have
$$H(r)(X,X)=I_{\gamma}(J^{\bot},J^{\bot})\le
I_{\gamma}(W^{\bot},W^{\bot})$$
$$=\frac{g_T(X^{\bot},X^{\bot})}{{\rm s}_c(r)^2}\int_0^{r}\left\{{\rm s}_c'(t)^2
-K(T(t);W(t)){\rm s}_c(t)^2\right\}dt$$
$$\le\frac{g_T(X^{\bot},X^{\bot})}{{\rm s}_c(r)^2}\int_0^{r}\left\{{\rm s}_c'(t)^2
-c{\rm s}_c(t)^2\right\}dt={\rm ct}_c(r)g_T(X^{\bot},X^{\bot}),$$
so we are done.\\

\Section{The Laplacian Comparison Theorems}{The Laplacian Comparison Theorems}
In this section we shall derive some Laplacian
comparison theorems for distance function. First of all, by Lemma
3.3 and Theorem 4.1 we have \\
{\bf Theorem 5.1} {\sl Let $(M,F,d\mu)$ be a Finsler $n$-manifold,
$r=d_F(p,\cdot)$, the distance function from a fixed point $p$.
Suppose that the flag curvature of $M$ satisfies $K(V;W)\le c$ for
any $V,W\in TM$. Then the following holds whenever $r$ is smooth.}
$$\triangle r\ge(n-1){\rm ct}_c(r)-\|{\bf S}\|,$$ {\sl where $\|{\bf
S}\|$ is the pointwise norm function of S-curvature which is
defined by}
$$\|{\bf S}\|_x=\sup_{X\in T_xM\backslash 0}\frac{{\bf
S}(X)}{F(X)}.$$ When $M$ has nonpositive flag curvature we have
the following Laplacian comparison theorem in terms of Ricci
curvature.\\
{\bf Theorem 5.2} {\sl Let $(M,F,d\mu)$ be a Finsler $n$-manifold
with nonpositive flag curvature. If the Ricci curvature of $M$
satisfies $Ric_M\le c<0$, then the following holds whenever $r$ is
smooth.}
$$\triangle r\ge {\rm ct}_c(r)-\|{\bf S}\|.$$
{\bf Proof.} We need only to prove tr$_{\nabla r}H(r)\ge {\rm
ct}_c(r)$. Suppose that $r$ is smooth at $q\in M$, then $r$ is
also smooth near $q$. Let $S_p(r(q))$ be the forward geodesic
sphere of radius $r(q)$ centered at $p$. We choose the local
$g_T$-orthonormal frame $E_1,\cdots,E_{n-1}$ of $S_p(r(q))$ near
$q$, here $T=\nabla r$. By parallel transformation along geodesic
rays we get local vector fields $E_1,\cdots,E_{n-1},E_n=T$ of $M$.
Then for any $1\le i,j\le n-1$, we have by (2.1)-(2.3) and (3.2),
$$\frac{d}{dr}\left(
H(r)(E_i,E_j)\right)=\frac{d}{dr}g_T\left(\nabla^T_{E_i}T,E_j\right)
=g_T\left(\nabla^T_T\nabla^T_{E_i}T,E_j\right)$$
$$=g_T\left(R^T(T,E_i)T,E_j\right)+g_T\left(\nabla^T_{[T,E_i]}T,E_j\right)$$
$$=-g_T\left(R^T(E_i,T)T,E_j\right)-g_T\left(\nabla^T_{\nabla^T_{E_i}T}T,E_j\right)$$
$$=-g_T\left(R^T(E_i,T)T,E_j\right)-\sum_kg_T\left(\nabla^T_{E_i}T,E_k\right)
g_T\left(\nabla^T_{E_k}T,E_j\right),$$ and consequently,
$$\frac{d}{dr}{\rm tr}_{\nabla r}H(r)=-Ric(\nabla
r)-\sum_{i,j}\left(H(r)(E_i,E_j)\right)^2.\eqno{(5.1)}$$ Since $M$
has nonpositive flag curvature, it is easy to see from Theorem 4.1
that the eigenvalues of $H(r)$ are nonnegative, which implies that
$$\sum_{i,j}\left(H(r)(E_i,E_j)\right)^2\le\left({\rm tr}_{\nabla
r}H(r)\right)^2,$$ and (5.1) can be rewritten as
$$\frac{d}{dr}{\rm tr}_{\nabla r}H(r)\ge-c-\left({\rm tr}_{\nabla
r}H(r)\right)^2.\eqno{(5.2)}$$ Note that in this case $c<0$, and
${\rm ct}_c(r)=\sqrt{-c}\cdot {\rm cotanh}(\sqrt{-c}r)$, from
(5.2) we have $$\frac{d}{dr}\left({\rm tr}_{\nabla r}H(r)-{\rm
ct}_c(r)\right)\ge-\left({\rm tr}_{\nabla r}H(r)\right)^2+{\rm
ct}_c(r)^2.\eqno{(5.3)}$$ Putting
$$A={\rm tr}_{\nabla r}H(r)-{\rm ct}_c(r),
B={\rm tr}_{\nabla r}H(r)+{\rm ct}_c(r),$$ then (5.3) becomes
$$\frac{dA}{dr}+AB\ge 0.\eqno{(5.4)}$$ We have again by the
nonpositivity of flag curvature and Theorem 4.1 that
$${\rm tr}_{\nabla r}H(r)\ge\frac{n-1}{r},$$ which implies that
there exist small $\varepsilon>0$ so that
$$A(r)\ge\frac{n-1}{r}-{\rm ct}_c(r)\ge 0,\quad\forall r\in
(0,\varepsilon].\eqno{(5.5)}$$ On the other hand, from (5.4) we
have
$$\frac{d}{dr}\left(A(r)\exp\left(\int_{\varepsilon}^rB(\tau)d\tau\right)
\right)\ge 0, $$ which yields
$$A(r)\exp\left(\int_{\varepsilon}^rB(\tau)d\tau\right)\ge
A(\varepsilon)\ge 0,$$ so we are done.\\
{\bf Remark} In the Riemannian case, ${\bf S}=0$, and Theorem 5.2
was obtained by [Ding] (see also [Xin1]).\\
\hspace*{0.7cm} For the case where the curvature is bounded from
below, we have the following comparison theorem.\\
{\bf Theorem 5.3} {\sl Let $(M,F,d\mu)$ be a Finsler $n$-manifold
with Ricci curvature satisfying $Ric_M\ge (n-1)c$. Then the
following holds whenever $r$ is smooth.}
$$\triangle r\le (n-1){\rm ct}_c(r)+\|{\bf S}\|. \eqno{(5.6)}$$
{\bf Proof.} Let $r=d_F(p,\cdot)$ is smooth at $q\in M$, and
$\gamma:[0,r(q)]\rightarrow M$ be the unit-speed geodesic from $p$
to $q$, and $T=\dot{\gamma}$. Let $e_1,\cdots,e_{n-1},e_n=T$ be
the $g_T$-orthonormal basis of $T_qM$. By parallel transformation
along $\gamma$ we obtain the parallel vector fields
$E_1(t),\cdots, E_n(t)$ along $\gamma$. For $1\le i\le n-1$, let
$J_i$ be the unique Jacobi field along $\gamma$ such that
$J_i(0)=0, J_i(r(q))=e_i$, and $W_i(t)=\frac{{\rm s}_c(t)}{{\rm
s}_c(r(q))}E_i(t)$, where ${\rm s}_c(t)$ is defined by (4.2).
Clearly, we have $W_i(0)=J_i(0)=0,W_i(r(q))=J_i(r(q))$. Thus from
(4.1) and the basic index lemma (see [BCS],page 182) we have
$${\rm tr}_{\nabla
r}(H(r))\mid_q=\sum_{i=1}^nH(r)(e_i,e_i)=\sum_{i=1}^{n-1}I_{\gamma}(J_i,J_i)
\le\sum_{i=1}^{n-1}I_{\gamma}(W_i,W_i)$$
$$=\frac{1}{{\rm s}_c(r(q))^2}\int_0^{r(q)}\left\{(n-1){\rm s}_c'(t)^2
-Ric(T(t)){\rm s}_c(t)^2\right\}dt$$
$$\le\frac{1}{{\rm s}_c(r(q))^2}\int_0^{r(q)}\left\{(n-1){\rm s}_c'(t)^2
-(n-1)c{\rm s}_c(t)^2\right\}dt=(n-1){\rm ct}_c(r(q)),$$ which
together with Lemma 3.3 yields (5.6).\\

\Section{Volume Comparison Theorems}{Volume Comparison Theorems}
In this section we shall use the Laplacian
comparison theorems to derive some volume comparison theorems for
Finsler manifolds.\\
\hspace*{0.7cm} Let $(M,F,d\mu)$ be a Finsler $n$-manifold. Fix
$p\in M$, let $I_p=\{ v\in T_pM: F(v)=1\}$ be the indicatrix at
$p$. For $v\in I_p$, the {\sl cut-value} $c(v)$ is defined by
$$c(v):=\sup\{t>0:d_F(p,\exp_p(tv))=t\}.$$ Then, we can define the
{\sl tangential cut locus} ${\bf C}(p)$ of $p$ by ${\bf C}(p):=\{
c(v)v:c(v)<\infty,v\in I_p\}$, the {\sl cut locus} $C(p)$ of $p$
by $C(p)=\exp_p{\bf C}(p)$, and the injectivity radius $i_p$ at
$p$ by $i_p=\inf\{c(v):v\in I_p\}$, respectively. It is known that
$C(p)$ has zero Hausdorff measure in $M$. Also, we set ${\bf
D}_p=\{tv:0\le t<c(v),v\in I_p\}$ and $D_p=\exp_p{\bf D}_p$. It is
known that ${\bf D}_p$ is the largest domain, starlike with
respect to the origin of $T_pM$, for which $\exp_p$ restricted to
that domain
is a diffeomorphism, and $D_p=M\backslash C(p)$.\\
\hspace*{0.7cm} Let $B_p(R)$ be the forward geodesic ball of $M$
with radius $R$ centered at $p$. The volume of $B_p(R)$ with
respect to $d\mu$ is defined by
$${\rm vol}(B_p(R))=\int_{B_p(R)}d\mu.$$
In order to compute the volume, we need the polar coordinates on
$D_p$. Let $\theta^{\alpha},\alpha=1,\cdots,n-1$ be the local
coordinates that are intrinsic to $I_p$. For any $q\in D_p$, the
polar coordinates of $q$ is defined by
$(r,\theta)=(r(q),\theta^1(q),\cdots,\theta^{n-1}(q))$, where
$r(q)=F(v),\theta^{\alpha}(q)=\theta^{\alpha}(\frac{v}{F(v)})$,
and $v=\exp_p^{-1}(q)$. Then by the Gauss lemma (see [BCS], page
140), the unit radial coordinate vector $\frac{\partial}{\partial
r}$ is $g_{\frac{\partial}{\partial r}}$-orthogonal to coordinate
vectors $\frac{\partial}{\partial\theta^{\alpha}}$ for
$\alpha=1,\cdots, n-1$. Therefore, writing
$d\mu=\sigma(r,\theta)dr\wedge
d\theta^1\wedge\cdots\wedge\theta^{n-1}:=\sigma(r,\theta)dr\wedge
d\theta$, we have, from (3.8),
$$\Delta r=\frac{\partial}{\partial r}\log\sigma.\eqno{(6.1)}$$
For $r>0$, let {\bf D}$_p(r)\subset I_p$ be defined by
$${\bf D}_p(r)=\{v\in I_p:rv\in {\bf D}_p\}.$$ It is easy to know
that ${\bf D}_p(r_1)\subset{\bf D}_p(r_2)$ for $r_1>r_2$ and
${\bf D}_p(r)=I_p$ for $r<i_p$. Since $C(p)$ has zero Hausdorff
measure in $M$, we have
$${\rm vol}(B_p(R))=\int_{B_p(R)}d\mu=\int_{{B_p(R)}\cap
D_p}d\mu$$
$$=\int_{\exp_p^{-1}(B_p(R))\cap{\bf
D}_p}\exp_p^*(d\mu)=\int_0^Rdr\int_{{\bf
D}_p(r)}\sigma(r,\theta)d\theta.\eqno{(6.2)}$$ For real numbers
$c,\Lambda$ and positive integer $n$, let
$$V_{c,\Lambda,n}(r)={\rm vol}(S^{n-1}(1))\int_0^re^{\Lambda t}{\rm s}_c(t)^{n-1}dt.\eqno{(6.3)}$$
We have\\
{\bf Theorem 6.1} {\sl Let $(M,F,d\mu)$ be a complete Finsler
$n$-manifold which satisfies $K(V;W)\le c$ and $\|{\bf
S}\|\le\Lambda$. Then the function}
$$\frac{{\rm vol}(B_p(r))}{V_{c,-\Lambda,n}(r)}$$ {\sl is monotone
increasing for $0<r\le i_p$, where $i_p$ is the injectivity radius
of $p$. In particular, for $d\mu=dV_F$, the Busemann-Hausdorff
volume form, one has }
$${\rm vol}(B_p(r)))\ge V_{c,-\Lambda,n}(r),\quad r\le i_p.\eqno{(6.4)}$$
{\bf Proof.} By (6.1), Theorem 5.1 and the assumptions of the
theorem, we have
$$\frac{\partial}{\partial r}\log\sigma\ge(n-1){\rm
ct}_c(r)-\Lambda=\frac{d}{dr}\log\left( e^{-\Lambda r}{\rm
s}_c(r)^{n-1}\right),$$ namely, the function
$$\frac{\sigma(r,\theta)}{e^{-\Lambda r}{\rm
s}_c(r)^{n-1}}$$ is monotone increasing in $r$ for any $\theta$.
Let $$\sigma_p(r)=\int_{{\bf D}_p(r)}\sigma(r,\theta)d\theta,\quad
\sigma_{c,-\Lambda,n}(r)={\rm vol}(S^{n-1}(1))e^{-\Lambda r}{\rm
s}_c(r)^{n-1}.$$ Then from (6.2) and (6.3) we have
$${\rm vol}(B_p(r)))=\int_0^r\sigma_p(t)dt,\quad
V_{c,-\Lambda,n}(r)= \int_0^r\sigma_{c,-\Lambda,n}(t)dt.$$ Noting
that {\bf D}$_p(r)=I_p$ for $r<
i_p,\frac{\sigma_p(r)}{\sigma_{c,-\Lambda,n}(r)}$ is also monotone
increasing for $r\le
 i_p$. Thus by the standard argument [Ch2],
the function
$$\frac{\displaystyle\int_0^r\sigma_p(t)dt}{\displaystyle\int_0^r\sigma_{c,-\Lambda,n}(t)dt}=
\frac{{\rm vol}(B_p(r))}{V_{c,-\Lambda,n}(r)}$$ is still monotone
increasing for $r\le i_p$. From [Sh3] we see that for $d\mu=dV_F$,
$$\lim_{r\rightarrow 0}\frac{{\rm vol}(B_p(r))}{V_{c,-\Lambda,n}(r)}=1,$$ thus we have (6.4).\\
\hspace*{0.7cm} The following theorem can be shown similarly by
use of Theorem 5.2.\\
{\bf Theorem 6.2} {\sl Let $(M,F,d\mu)$ be a complete and simply
connected Finsler $n$-manifold with nonpositive flag curvature. If
the Ricci curvature of $M$ satisfies $Ric_M\le c<0$ and $\|{\bf
S}\|\le\Lambda$, then the function}
$$\frac{{\rm vol}(B_p(r))}{V_{c,-\Lambda,2}(r)}$$ {\sl is monotone
increasing. In particular, for $d\mu=dV_F$, }
$${\rm vol}(B_p(r)))\ge \frac{{\rm vol}(B^n(1))}{{\rm vol}(B^2(1))}V_{c,-\Lambda,2}(r).$$
\hspace*{0.7cm} The following theorem was first obtained in [Sh3],
and here we provide another proof by use of Laplacian comparison theorem.\\
{\bf Theorem 6.3}$^{[Sh3]}$ {\sl Let $(M,F,d\mu)$ be a complete
Finsler $n$-manifold. Suppose that }
$$ Ric_M\ge(n-1)c,\quad\|{\bf S}\|\le \Lambda.$$ Then the function
$$\frac{{\rm vol}(B_p(r))}{V_{c,\Lambda,n}(r)}$$ {\sl is monotone decreasing in $r$. In particular, for
$d\mu=dV_F$, }
$${\rm vol}(B_p(r))\le V_{c,\Lambda,n}(r).$$
{\bf Proof.} By (6.1), Theorem 5.3 and the assumptions of the
theorem we have
$$\frac{\partial}{\partial r}\log\sigma\le(n-1){\rm
ct}_c(r)+\Lambda=\frac{d}{dr}\log\left( e^{\Lambda r}{\rm
s}_c(r)^{n-1}\right),$$ thus the function
$$\frac{\sigma(r,\theta)}{e^{\Lambda r}{\rm
s}_c(r)^{n-1}}$$ is monotone decreasing.  Noting that ${\bf
D}_p(R)\subset{\bf D}_p(r)$ for $R>r>0$, we have for $R>r>0$,
$$\frac{\sigma_p(r)}{\sigma_{c,\Lambda,n}(r)}=\frac{1}{{\rm
vol}(S^{n-1}(1))} \int_{{\bf
D}_p(r)}\frac{\sigma(r,\theta)}{e^{\Lambda r}{\rm
s}_c(r)^{n-1}}d\theta\ge\frac{1}{{\rm vol}(S^{n-1}(1))} \int_{{\bf
D}_p(R)}\frac{\sigma(r,\theta)}{e^{\Lambda r}{\rm
s}_c(r)^{n-1}}d\theta$$
$$\ge\frac{1}{{\rm
vol}(S^{n-1}(1))} \int_{{\bf
D}_p(R)}\frac{\sigma(R,\theta)}{e^{\Lambda R}{\rm
s}_c(R)^{n-1}}d\theta=\frac{\sigma_p(R)}{\sigma_{c,\Lambda,n}(R)},$$
namely, $\frac{\sigma_p(r)}{\sigma_{c,\Lambda,n}(r)}$ is also
monotone decreasing. Now the theorem can be verified easily.

\Section{The First Eigenvalue}{The First Eigenvalue}
In this section we shall study the first eigenvalue for Finsler
manifolds and proof some Mckean type theorems. We need some lemmas.\\
{\bf Lemma 7.1} {\sl Let $(M,F)$ be a Finsler manifold with finite
reversibility $\lambda$, then $\mid\xi(X)\mid\le\lambda
F^*(\xi)F(X)$ for any} $X\in TM, \xi\in T^*M$.\\
{\bf Proof.} By the definition of $F^*$ we actually have
$\xi(X)\le F^*(\xi)F(X)$. On the other hand, from the definition
of reversibility one has $-\xi(X)=\xi(-X)\le
F^*(\xi)F(-X)\le\lambda
F^*(\xi)F(X)$, so the lemma follows.\\
\hspace*{0.7cm} Now let $(M,F,d\mu)$ be a Finsler $n$-manifold,
$\Omega\subset M$ a domain with compact closure and nonempty
boundary $\partial\Omega$. The first eigenvalue
$\lambda_1(\Omega)$ of $\Omega$ is defined by (see [Sh2], page
176)
$$\lambda_1(\Omega)=\inf_{f\in L^2_{1,0}(\Omega)
\backslash\{0\}}\left\{\frac{\displaystyle\int_{\Omega}\left(F^*(df)\right)^2d\mu}{\displaystyle\int_{\Omega}f^2d\mu}\right\},$$
where $L^2_{1,0}(\Omega) $ is the completion of $C^{\infty}_0$
with respect to the norm
$$\parallel\varphi\parallel^2_{\Omega}=\int_{\Omega}\varphi^2d\mu+\int_{\Omega}\left(F^*(df)\right)^2d\mu.$$
If $\Omega_1\subset\Omega_2$ are bounded domains, then
$\lambda_1(\Omega_1)\ge\lambda_1(\Omega_2)\ge 0$. Thus, if
$\Omega_1\subset\Omega_2\subset\cdots\subset M$ be bounded domains
so that $\bigcup\Omega_i=M$, then the following limit
$$\lambda_1(M)=\lim_{i\rightarrow\infty}\lambda_1(\Omega_i)\ge 0$$
exists, and it is independent of the choice of $\{\Omega_i\}$. we
have the following lemma which is crucial in this section.\\
{\bf Lemma 7.2} {\sl Let $(M,F,d\mu)$ be a Finsler manifold with
finite reversibility $\lambda$, $\Omega\subset M$ a domain with
compact closure and nonempty boundary, and $X$ a vector field on
$\Omega$ so that $\|X\|_{\infty}=\sup_{\Omega}F(X)<\infty$ and
$\inf_{\Omega}{\rm div}(X)>0$. Then}
$$\lambda_1(\Omega)\ge\left[\frac{\inf_{\Omega}{\rm
div}X}{2\lambda\|X\|_{\infty}}\right]^2,\eqno{(7.1)}$$ {\bf
Proof.} Let $f\in C_0^{\infty}$, then vector field $f^2X$ has
compact support in $\Omega$. Now computing the divergence of
$f^2X$ we have by Lemma 7.1,
$${\rm div}(f^2X)=2fX(f)+f^2{\rm div}X$$
$$\ge-2\lambda\mid f\mid\cdot\sup_{\Omega}F(X)\cdot F^*(df)+\inf_{\Omega}{\rm
div}X\cdot f^2.\eqno{(7.2)}$$ Using the inequality
$$-2\mid f\mid\cdot F^*(df)\ge-\varepsilon
f^2-\frac{1}{\varepsilon} \left(F^*(df)\right)^2$$ for all
$\varepsilon>0$, we have from (7.2) that
$${\rm div}(f^2X)\ge\lambda\cdot\sup_{\Omega}F(X)\cdot\left(-\varepsilon
f^2-\frac{1}{\varepsilon}
\left(F^*(df)\right)^2\right)+\inf_{\Omega}{\rm div}X\cdot
f^2.\eqno{(7.3)}$$ Integrating (7.3) on $\Omega$ and using Lemma
3.2 we have
$$0=\int_{\Omega}{\rm div}(f^2X)d\mu$$
$$\ge\lambda\cdot\|X\|_{\infty}\int_{\Omega}\left(-\varepsilon
f^2-\frac{1}{\varepsilon}
\left(F^*(df)\right)^2\right)d\mu+\inf_{\Omega}{\rm
div}X\cdot\int_{\Omega}f^2d\mu.$$ Therefore,
$$\int_{\Omega}\left(F^*(df)\right)^2d\mu\ge\frac{\varepsilon}{\lambda\|X\|_{\infty}}
\left(\inf_{\Omega}{\rm
div}X-\lambda\cdot\|X\|_{\infty}\cdot\varepsilon\right)\int_{\Omega}
f^2d\mu.$$ Choosing $\varepsilon=\inf_{\Omega}{\rm
div}X/(2\lambda\cdot\|X\|_{\infty})$ we have
$$\int_{\Omega}\left(F^*(df)\right)^2d\mu\ge\left[\frac{\inf_{\Omega}{\rm
div}X}{2\lambda\|X\|_{\infty}}\right]^2\int_{\Omega}
f^2d\mu.\eqno{(7.4)}$$ Since (7.4) holds for any $f\in
C^{\infty}_0(\Omega)$, we have (7.1).\\
\hspace*{0.7cm} Now we are in the position to prove the first main
result of this section.\\
{\bf Theorem 7.3} {\sl Let $(M,F,d\mu)$ be a Finsler $n$-manifold
with finite reversibility $\lambda$ and flag curvature $K(V;W)\le
c$ for any $V,W\in TM$. Let $B_p(R)$ be the forward geodesic ball
of $M$ with radius $R$ centered at $p$, and $R<i_p$, where $i_p$
denotes the injectivity radius about $p$. Suppose that $(n-1){\rm
ct}_c(R)-\sup_{B_p(R)}\|{\bf S}\|>0$, then }
$$\lambda_1(B_p(R))\ge \left[\frac{(n-1){\rm ct}_c(R)-\sup_{B_p(R)}\|{\bf
S}\|}{2\lambda}\right]^2.\eqno{(7.5)}$$ {\bf Proof.} For
$R>\varepsilon>0$, let
$\Omega_{\varepsilon}=B_p(R)\backslash\overline{B_p(\varepsilon)}$.
Then $r=d_F(p,\cdot)$ is smooth on $\Omega_{\varepsilon}$, and
thus $X=\nabla r$ is a smooth vector field on
$\Omega_{\varepsilon}$. Noting that $F(X)=F(\nabla r)=1$ and
div$X=\triangle r$, we deduce from Theorem 5.1 and Lemma 7.2 that
$$\lambda_1(\Omega_{\varepsilon})\ge \left[\frac{(n-1){\rm ct}_c(R)-\sup_{B_p(R)}\|{\bf
S}\|}{2\lambda}\right]^2.$$ Letting $\varepsilon\rightarrow 0$ we
get (7.5).\\
\hspace*{0.7cm} By Theorem 7.3 we have the following result which
is the Finsler version of Mckean's theorem.\\
{\bf Theorem 7.4} {\sl Let $(M,F,d\mu)$ be a complete noncompact
and simply connected Finsler $n$-manifold with finite
reversibility $\lambda$ and flag curvature $K(V;W)\le-a^2(a>0)$.
If $\sup_M\|{\bf S}\|<(n-1)a$, then }
$$\lambda_1(M)\ge\frac{\left((n-1)a-\sup_M\|{\bf S}\|\right)^2}{4\lambda^2}.$$
\hspace*{0.7cm} The following result can be verified by use of
Theorem 5.2 which is another Finsler version of Mckean's theorem
in term of Ricci curvature.\\
{\bf Theorem 7.5} {\sl Let $(M,F,d\mu)$ be a complete noncompact
and simply connected Finsler $n$-manifold with finite
reversibility $\lambda$ and nonpositive flag curvature. If
$Ric_M\le-a^2(a>0)$ and $\sup_M\|{\bf S}\|<a$, then}
$$\lambda_1(M)\ge\frac{\left(a-\sup_M\|{\bf S}\|\right)^2}{4\lambda^2}.$$

\Section{On Curvature and Fundamental Group}{On Curvature and Fundamental Group}
In 1968 Milnor [Mi] studied the curvature and
fundamental group of Riemannian manifolds and proved that the
fundamental group of a compact Riemannian manifold with strictly
negative sectional curvature has at least exponential growth. The key in the
proof is that the fundamental group can be identified with the
deck transformation group of the universal covering space, and any
geodesic ball in universal covering space can be covered by the union of a number of
translations of the fundamental domain. Combining with the
estimate of the volume growth Milnor was able to obtain his
result. His result was generalized in [Y] and [Xin2]. Milnor's idea can be generalized
to the Finsler setting. As the first step,
we have, from Theorems 6.1 and 6.2,\\
{\bf Lemma 8.1} {\sl Let $(M,F,d\mu)$ be a simply connected and
complete Finsler $n$-manifold with $\|{\bf S}\|\le\Lambda$.
Suppose that one of the following two conditions holds: } \\
(i) {\sl the flag curvature of $M$ satisfies $K(V;W)\le-a^2$ with
$a>\Lambda/(n-1)$;}\\
(ii) {\sl $M$ has nonpositive flag curvature and $Ric_M\le-a^2$
with $a>\Lambda$.}\\
{\sl Then the volume of the forward geodesic ball of $M$ grows at
least exponentially.}\\
\hspace*{0.7cm} For the universal covering space of a Finsler
manifold, we can endow the covering space with a pulled-back Finsler
metric so that the covering map is a local isometry. With Lemma 8.1
at hand, we can prove the following theorem by the almost same
argument as in [Mi].\\
{\bf Theorem 8.2} {\sl Let $(M,F,d\mu)$ be a compact Finsler
$n$-manifold with $\|{\bf S}\|\le\Lambda$.
Suppose that one of the following two conditions holds: } \\
(i) {\sl the flag curvature of $M$ satisfies $K(V;W)\le-a^2$ with
$a>\Lambda/(n-1)$;}\\
(ii) {\sl $M$ has nonpositive flag curvature and $Ric_M\le-a^2$
with $a>\Lambda$.}\\
{\sl Then the fundamental group of $M$ grows at least
exponentially.}\\


\begin{thebibliography}{99}
\bibitem[AL]{1} Antonelli P.L. and Lackey B., eds, Proc. Conf. on
Finsler Laplacians, Kluwer Academic Press, Netherlands, 1998.
\bibitem[BCS]{2} Bao D., Chern S.S. and Shen Z., An introduction to Riemannian-Finsler geometry,
GTM 200, Springer-Verlag, 2000.
\bibitem[Ce]{3} Centore P., Finsler Laplacians and minimal-energy
maps, Int. J. of Math., 11(2000),1-13.
\bibitem[Ch1]{4} Chern S.S., Local equivlence and Euclidean
connections in Finsler spaces, Sci. Rep. Nat. Tsing Hua Univ. Ser.
A5(1948),95-121; or Selected Papers, II, 194-212, Springer 1989.
\bibitem[Ch2]{5} Chavel I., Riemannian geometry, a modern
introduction, Camb. Univ. Press, 1993.
\bibitem[Ding]{6} Ding Q., A new Laplacian comparison theorem and
the estimate of eigenvalues, Chin. Ann. of Math., 15B(1994),35-42.
\bibitem[Fin]{7} Finsler P., \"{U}ber Kurven und Fl\"{a}chen
in allgemeinen R\"{a}umen, Dissertation, G\"{o}ttingen 1918.
\bibitem[Mc]{8} Mckean H.P., An upper bound fot the spectrum of
$\triangle$ on a manifold of negative curvature, J. Diff. Geom.,
4(1970),359-366.
\bibitem[Mi]{9} Milnor J., A note on curvature and fundamental
group, J. Diff. Geom., 2(1968), 1-7.
\bibitem[Ra]{10} Rademacher H.B., A sphere theorem for
non-reversible Finsler metrics, Math. Ann., 328(2004), 373-387.
\bibitem[Sh1]{11} Shen Z., On Finsler geometry of submanifolds, Math. Ann., 311(1998), 549-576.
\bibitem[Sh2]{12} Shen Z., Lectures on Finsler geometry, World Sci.,
2001, Singapore.
\bibitem[Sh3]{13} Shen Z., Volume comparison and its applications
in Riemann-Finsler geometry, Adv. in Math., 128(1997), 306-328.
\bibitem[Xin1]{14} Xin Y.L., Geometry of harmonic maps,
Birkh\"{a}user PNLDE 23,1996.
\bibitem[Xin2]{15} Xin Y. L., Ricci curvature and fundamental group, Preprint.
\bibitem[Y]{16} Yang, Yi-Hu, On the growth of fundamental groups
on nonpositive curvature manifolds, Bull. Austral. Math. 54 (1996), 483-487.
\end{thebibliography}
\end{document}